\documentclass[12pt,twoside]{article}

\begin{document}

\centerline{}

\centerline {\Large{\bf Gateaux and Fr$\acute{e}$chet Derivative in Intuitionistic }}
\centerline{}
\centerline{\Large{\bf  Fuzzy Normed Linear spaces}}

\centerline{}
\newcommand{\mvec}[1]{\mbox{\bfseries\itshape #1}}

\centerline{\bf {B. Dinda, T.K. Samanta and U.K. Bera}}

\centerline{}

\centerline{Department of Mathematics,}
\centerline{Mahishamuri Ramkrishna
Vidyapith, West Bengal, India. }
\centerline{e-mail: bvsdinda@gmail.com}

\centerline{Department of Mathematics, Uluberia
College, West Bengal, India.}
\centerline{e-mail: mumpu$_{-}$tapas5@yahoo.co.in}

\centerline{Department of Mathematics, City
College, Kolkata, India-700009.}
\centerline{e-mail: uttamkbera@gmail.com}
\centerline{}

\newtheorem{Theorem}{\quad Theorem}[section]

\newtheorem{definition}[Theorem]{\quad Definition}

\newtheorem{theorem}[Theorem]{\quad Theorem}

\newtheorem{remark}[Theorem]{\quad Remark}

\newtheorem{corollary}[Theorem]{\quad Corollary}

\newtheorem{note}[Theorem]{\quad Note}

\newtheorem{lemma}[Theorem]{\quad Lemma}

\newtheorem{example}[Theorem]{\quad Example}

\begin{abstract}
\textbf{\emph{Intuitionistic Fuzzy derivative, Intuitionistic Fuzzy Gateaux derivative, Intuitionistic Fuzzy Fr\'{e}chet derivative are defined and a few of their properties are studied. The relation between Intuitionistic Fuzzy Gateaux derivative and Intuitionistic Fuzzy Fr\'{e}chet derivative are emphasized.}}
\end{abstract}
${\bf Keywords:}$  \emph{Intuitionistic fuzzy differentiation, intuitionistic fuzzy continuity, intuitionistic fuzzy Gateaux derivative, intuitionistic fuzzy Fr$\acute{e}$chet derivative.}

\section{Introduction}
Fuzzy set theory is a useful tool to describe the situation in which data are imprecise or vague or uncertain. Intuitionistic fuzzy set theory handle the situation by attributing a degree of membership and a degree of non-membership to which a certain object belongs to a set. It has a wide range of application in the field of population dynamics \cite{Barros LC}, chaos control \cite{Fradkov}, computer programming \cite{Giles}, medicine \cite{Barro} etc. \\
The concept of intuitionistic fuzzy set, as a generalisation of fuzzy sets \cite{zadeh}
was introduced by Atanassov in \cite{Atanassov}. The concept of fuzzy norm was introduced by Katsaras
\cite{Katsaras} in 1984. In 1992, Felbin\cite{Felbin1} introduced
the idea of fuzzy norm on a linear space. Cheng-Moderson
\cite{Shih-chuan} introduced another idea of fuzzy norm on a linear
space whose associated metric is same as the associated metric of
Kramosil-Michalek \cite{Kramosil}. Latter on Bag and Samanta
\cite{Bag} modified the definition of fuzzy norm of
Cheng-Moderson \cite{Shih-chuan} and established the concept of
continuity and boundednes of a linear operator with respect to
their fuzzy norm in \cite{Bag1}.\\
Many authors in \cite{Ernesto,Ferrarro,Kuhn,Puri} discuss fuzzy derivatives in many approach. After studying continuities and boundedness  of linear operator in fuzzy environment in \cite{dinda,Samanta1,dinda1,dinda2,dinda3}, we introduce intuitionistic fuzzy Gateaux derivative and intuitionistic fuzzy Fr$\acute{e}$chet derivative of linear operator.\\
In this paper we define intuitionistic fuzzy derivative in $\mathbf{R}$, intuitionistic fuzzy Gateaux derivative and intuitionistic fuzzy Fr$\acute{e}$chet derivative of linear operator and we study some of their properties. Thereafter we show that in $\,\mathbf{R}\,$ intuitionistic fuzzy derivative and intuitionistic fuzzy Fr$\acute{e}$chet derivative are equivalent. We also show that intuitionistic fuzzy Fr$\acute{e}$chet derivative implies intuitionistic fuzzy Gateaux derivative.

\section{Preliminaries}
We quote some definitions and statements of a few theorems
which will be needed in the sequel.

\begin{definition}
\cite{Schweizer} A binary operation \, $\ast \; : \; [\,0 \; , \;
1\,] \; \times \; [\,0 \; , \; 1\,] \;\, \longrightarrow \;\, [\,0
\; , \; 1\,]$ \, is continuous \, $t$ - norm if \,$\ast$\, satisfies
the
following conditions \, $:$ \\
$(\,i\,)$ \hspace{0.5cm} $\ast$ \, is commutative and associative ,
\\ $(\,ii\,)$ \hspace{0.4cm} $\ast$ \, is continuous , \\
$(\,iii\,)$ \hspace{0.2cm} $a \;\ast\;1 \;\,=\;\, a \hspace{1.2cm}
\forall \;\; a \;\; \varepsilon \;\; [\,0 \;,\; 1\,]$ , \\
$(\,iv\,)$ \hspace{0.2cm} $a \;\ast\; b \;\, \leq \;\, c \;\ast\; d$
\, whenever \, $a \;\leq\; c$  ,  $b \;\leq\; d$  and  $a \,
, \, b \, , \, c \, , \, d \;\, \varepsilon \;\;[\,0 \;,\; 1\,]$.
\end{definition}
A few examples of continuous t-norm are $\,a\,\ast\,b\,=\,ab,\;\,a\,\ast\,b\,=\,min\{a,b\},\;\,a\,\ast\,b\,=\,max\{a+b-1,0\}$.

\begin{definition}
\cite{Schweizer}. A binary operation \, $\diamond \; : \; [\,0 \; ,
\; 1\,] \; \times \; [\,0 \; , \; 1\,] \;\, \longrightarrow \;\,
[\,0 \; , \; 1\,]$ \, is continuous \, $t$-conorm if
\,$\diamond$\, satisfies the
following conditions \, $:$ \\
$(\,i\,)\;\;$ \hspace{0.1cm} $\diamond$ \, is commutative and
associative ,
\\ $(\,ii\,)\;$ \hspace{0.1cm} $\diamond$ \, is continuous , \\
$(\,iii\,)$ \hspace{0.1cm} $a \;\diamond\;0 \;\,=\;\, a
\hspace{1.2cm}
\forall \;\; a \;\; \in\;\; [\,0 \;,\; 1\,]$ , \\
$(\,iv\,)$ \hspace{0.1cm} $a \;\diamond\; b \;\, \leq \;\, c
\;\diamond\; d$ \, whenever \, $a \;\leq\; c$  ,  $b \;\leq\; d$
 and  $a \, , \, b \, , \, c \, , \, d \;\; \in\;\;[\,0
\;,\; 1\,].$
\end{definition}
A few examples of continuous t-conorm are $\,a\,\ast\,b\,=\,a+b-ab,\;\,a\,\ast\,b\,=\,max\{a,b\},\;\,a\,\ast\,b\,=\,min\{a+b,1\}$.

\begin{definition}
\cite{Samanta} Let \,$\ast$\, be a continuous \,$t$-norm ,
\,$\diamond$\, be a continuous \,$t$-conorm  and \,$V$\, be a
linear space over the field \,$F \;(\, = \; \mathbf{R} \;\, or \;\,
\mathbf{C} \;)$. An \textbf{intuitionistic fuzzy norm} on \,$V$\,
is an object of the form  $\\A\;=\{\;(\,(\,x \;,\; t\,)
\;,\;\mu\,(\,x \;,\; t\,) \;,
\; \nu\,(\,x \;,\; t\,)\;)\;:
\;(\,x \;,\; t\,) \;\,\in\; V \;\times\;
\mathbf{R^{\,+}} \;\}$ , where $\mu \,,\, \nu\;are\; fuzzy\; sets
\;on \,$V$  \;\times\; \mathbf{R^{\,+}}$ , \,$\mu$\, denotes the
degree of membership and \,$\nu$\, denotes the degree of non -
membership \,$(\,x \;,\; t\,) \;\,\in\;\, V \;\times\;
\mathbf{R^{\,+}}$\, satisfying the following conditions $:$ \\\\
$(\,i\,)$ \hspace{0.10cm}  $\mu\,(\,x \;,\; t\,) \;+\; \nu\,(\,x
\;,\; t\,) \;\,\leq\;\, 1 \hspace{1.2cm} \forall \;\; (\,x \;,\;
t\,)
\;\,\in\;\, V \;\times\; \mathbf{R^{\,+}}\, ;$ \\
$(\,ii\,)$ \hspace{0.10cm}$\mu\,(\,x \;,\; t\,) \;\,>\;\, 0 \, ;$ \\
$(\,iii\,)$ $\mu\,(\,x \;,\; t\,) \;\,=\;\, 1\;$ if
and only if \, $x \;=\; \theta \,$, $\theta$ is null vector ; \\
$(\,iv\,)$\hspace{0.05cm} $\mu\,(\,c\,x \;,\; t\,) \;\,=\;\,
\mu\,(\,x \;,\; \frac{t}{|\,c\,|}\,)\;\;\;\;\forall\; c
\;\,\in\;\, F \, $ and $c \;\neq\; 0 \;;$ \\ $(\,v\,)$ \hspace{0.10cm} $\mu\,(\,x
\;,\; s\,) \;\ast\; \mu\,(\,y \;,\; t\,) \;\,\leq\;\, \mu\,(\,x
\;+\; y \;,\; s \;+\; t\,) \, ;$ \\ $(\,vi\,)$ \hspace{0.05cm}
$\mu\,(\,x \;,\; \cdot\,)$ is non-decreasing function of
$\;\mathbf{R^{\,+}}\,$ and  $\,\mathop {\lim }\limits_{t\;\, \to
\,\;\infty } \;\,\,\mu\,\left( {\;x\;,\;t\,} \right)=1 ;$
\\ $(\,vii\,)$ \hspace{0.10cm}$\nu\,(\,x \;,\; t\,) \;\,<\;\, 1 \, ;$ \\
$(\,viii\,)$ $\nu\,(\,x \;,\; t\,) \;\,=\;\, 0\;$  if
and only if  $\,x \;=\; \theta \, ;$ \\ $(\,ix\,)$
\hspace{0.05cm} $\nu\,(\,c\,x \;,\; t\,) \;\,=\;\, \nu\,(\,x \;,\;
\frac{t}{|\,c\,|}\,)\;\;\;\;\forall\; c
\;\,\in\;\, F \, $ and $c \;\neq\; 0 \;;$ \\ $(\,x\,)$ \hspace{0.15cm} $\nu\,(\,x
\;,\; s\,) \;\diamond\; \nu\,(\,y \;,\; t\,) \;\,\geq\;\, \nu\,(\,x
\;+\; y \;,\; s \;+\; t\,) \, ;$ \\ $(\,xi\,)$ \hspace{0.04cm}
$\nu\,(\,x \;,\; \cdot\,)$ is non-increasing function of \,
$\mathbf{R^{\,+}}\;$  and  $\;\mathop {\lim }\limits_{t\;\, \to
\,\;\infty } \;\,\,\nu\,\left( {\;x\;,\;t\,} \right)=0.$
\end{definition}

\begin{definition}
\cite{Samanta} If $A$ is an intuitionistic fuzzy norm on a linear
space $V$ then $(V\;,\;A)$ is called an intuitionistic fuzzy normed
linear space.
\end{definition}

 For the intuitionistic fuzzy normed linear space $\;(\,V \;,\; A\,)\;$,
 we further assume that $\;\mu,\, \nu,\, \ast,\, \diamond$\,
 satisfy the following axioms : \\
$(\,xii\,)\;\;\;\;\;\;\left. {{}_{a\;\; \ast \;\;a\;\; =
\;\;a}^{a\;\; \diamond \;\;a\;\; = \;\;a} \;\;}
\right\}\;\;\;\,\;\;\;\;\,,\;\;\,$ for all
$\;a\;\; \in \;\;[\,0\;\,,\;\,1\,].$ \\
$(\,xiii\,)\;\;\mu\,(\,x \;,\; t\,) \;>\; 0 \;\;\;\;,\;$ for all $
\;\; t \;>\; 0 \;\; \Rightarrow \;\; x \;=\;\theta\;.$ \\
$(\,xiv\,)\;\;\;\nu\,(\,x \;,\; t\,) \;<\; 1 \;\;\;\;\;\; ,\;$ for all
$\;\; t \;>\; 0 \;\; \Rightarrow \;\; x \;=\; \theta\;.$ \\
$(xv)\;$ For $\;x\;\neq\;\theta,\;\;\mu(x\,,\, .)$ is a continuous function of $\mathbf{R}$ and strictly increasing on the subset $\{\, t\;\,:\;\,0\;<\;\mu(x\,,\,t)\;<\;1 \, \}$ of $\mathbf{R}$.\\
$(xvi)\;$ For $\;x\;\neq\;\theta,\;\;\;\nu(x\,,\, .)$ is a continuous function
 of $\mathbf{R}$ and strictly decreasing on the
 subset $\{\,t\;\,:\;\,0\;<\;\nu(x\,,\,t)\;<\;1 \,\}$ of $\mathbf{R}$.

\begin{definition}
\cite{Samanta} A sequence $\{x_n\}_n$ in an intuitionistic fuzzy normed linear space $(V\,,\,A)$ is said to \textbf{converge} to $x\;\in\;V$ if for given $r>0,\;t>0,\;0<r<1$, there exist an integer $n_0\;\in\;\mathbf{N}$ such that \\
$\;\mu\,(\,x_n\,-\,x\,,\,t\,)\;>\;1\,-\,r\;\;$
 and $\;\;\nu\,(\,x_n\,-\,x\,,\,t\,)\;<\;r\;\;$  for all $n\;\geq \;n_0$.
\end{definition}

\begin{definition}
\cite{Samanta} Let, $(\;U\;,\;A\;)$ and $(\;V\;,\;B\;)$ be two
intuitionistic fuzzy normed linear space over the same field $F$. A
mapping $f$ from $(\;U\;,\;A\;)$ to $(\;V\;,\;B\;)$ is said to be \textbf{
intuitionistic fuzzy continuous} at $x_0\;\in\;U$, if for any given
$\epsilon\;>\;0\;,\alpha\;\in\;(0,1)\;,\exists\;\delta
\;=\delta(\alpha,\epsilon)\;>0\;,\beta\;=\beta(\alpha,\epsilon)\;\in\;(0,1)$
such that for all $x\;\in\;U$,
\[\mu_U(x-x_0 \;,\;\delta)\;>\;1-\beta\;\Rightarrow\;
\mu_V(f(x)-f(x_0) \;,\;\epsilon)\;>\;1-\alpha\] \[
\nu_U(x-x_0 \;,\;\delta)\;<\;\beta\;\Rightarrow\;
\nu_V(f(x)-f(x_0) \;,\;\epsilon)\;<\;\alpha\;. \]
\end{definition}

\section{Intuitionistic fuzzy Gateaux derivative}
In this section, we shall consider $\;(\,\mathbf{R},\mu_R,\nu_R,\ast,\diamond\,)\;$ as an intuitionistic fuzzy normed linear space over the field $\;\mathbf{R}\,$(the set of all real numbers).

\begin{definition}
 Let $\;(\,\mathbf{R},\mu_1,\nu_1,\ast,\diamond\,)\;$ and $\;(\,\mathbf{R},\mu_2,\nu_2,\ast,\diamond\,)\;$ be two
intuitionistic fuzzy normed linear space over the same field $\;\mathbf{R}.\;$ A
mapping $\;f\;$ from $\;(\,\mathbf{R},\mu_1,\nu_1,\ast,\diamond\,)\;$ to $\;(\,\mathbf{R},\mu_2,\nu_2,\ast,\diamond\,)\;$ is said to be {\bf
intuitionistic fuzzy differentiable} at $\;x_0\;\in\;\mathbf{R},\;$ if for any given
$\;\epsilon\;>\;0\;,\alpha\,\in\,(0,1)\;,\exists\;\delta
\,=\,\delta(\alpha,\epsilon)\;>0\;,\beta\;=\beta(\alpha,\epsilon)\;\in\;(0,1)$
such that for all $x(\neq\,x_0)\,\in\,\mathbf{R}$,
\[\mu_1(x-x_0 \;,\;\delta)\;>\;1-\beta\;\Rightarrow\;\mu_2\left(\,\frac{f(x)-f(x_0)}{x-x_0}\,-\,f^\prime (x_0) \,,\,\epsilon\,\right)\,>\,1-\alpha\]
 \[\nu_1(x-x_0 \;,\;\delta)\;<\;\beta\;\Rightarrow\;\nu_2\left(\,\frac{f(x)-f(x_0)}{x-x_0}\,-\,f ^\prime(x_0) \,,\,\epsilon\,\right)\,<\,\alpha\;. \]
 We denote intuitionistic fuzzy derivative of $\;f\;$ at $\;x_0\;$ by $f^\prime(x_0).$
\end{definition}

{\bf Alternative definition:} Let $\;(\,\mathbf{R},\mu_1,\nu_1,\ast,\diamond\,)\;$ and $\;(\,\mathbf{R},\mu_2,\nu_2,\ast,\diamond\,)\;$ be two
intuitionistic fuzzy normed linear space over the same field $\;\mathbf{R}.\;$ A
mapping $\;f\;$ from $\;(\,\mathbf{R},\mu_1,\nu_1,\ast,\diamond\,)\;$ to $\;(\,\mathbf{R},\mu_2,\nu_2,\ast,\diamond\,)\;$ is said to be {\bf
intuitionistic fuzzy differentiable} at $\;x_0\;\in\;\mathbf{R},\;$ if for every $\;t>0\;$
\[\mathop {\lim }\limits_{{\mu_1(x-x_0,t)}\;\, \to\,\;1 } \;\mu_2\left(\,\frac{f(x)-f(x_0)}{x-x_0}\,-\,f^\prime(x_0) \,,\,t\,\right)\,=\,1\]
\[\mathop {\lim }\limits_{{\nu_1(x-x_0,t)}\;\, \to\,\;0 } \;\nu_2\left(\,\frac{f(x)-f(x_0)}{x-x_0}\,-\,f^\prime(x_0) \,,\,t\,\right)\,=\,0\]
$f^\prime(x_0)\;$ is called intuitionistic fuzzy derivative of $\;f\;$ at $\;x_0\;.$

\begin{note}
It is easy to see that these two definitions are equivalent.
\end{note}

\begin{note}
If the intuitionistic fuzzy derivative of $\;f,\;$ be $\;f^\prime(x_0),\;$ the intuitionistic fuzzy derivative of $\;f^\prime(x_0)\;$ at $\;x_0\;$ is called second order intuitionistic fuzzy derivative of $\;f\;$ at $\;x_0\;$ and is denoted by $\;f^{\prime\prime}(x_0).\;$ Similarly, the n-th order intuitionistic fuzzy derivative of $\;f\;$ at $\;x_0\;$ exists if $\;{f}^{n-1}(x_0)\;$ is intuitionistic fuzzy differentiable at $\;x_0\;$ and this derivative is denoted by $\;{f}^{n}(x_0)\,.$
\end{note}

\begin{theorem}
Let $\;f\;:(\,\mathbf{R},\mu_1,\nu_1,\ast,\diamond\,)\;\rightarrow\;
(\,\mathbf{R},\mu_2,\nu_2,\ast,\diamond\,)\;$ and $\;g\;:(\,\mathbf{R},\mu_1,\nu_1,\ast,\diamond\,)\;\rightarrow\;
(\,\mathbf{R},\mu_2,\nu_2,\ast,\diamond\,)\;$ be intuitionistic fuzzy differentiable at $\;x_0,\;(\,\mathbf{R},\mu_1,\nu_1,\ast,\diamond\,)\;$
and $\;(\,\mathbf{R},\mu_2,\nu_2,\ast,\diamond\,)\;$ satisfies the condition $\;(xii).\;$ Then for $\;K\in\mathbf{R},\;\;Kf+g\;$ is intuitionistic fuzzy differentiable at $\;x_0\;$ and $\;(Kf+g)^\prime(x_0)\,=\,K{f}^\prime(x_0)\,+\,{g}^\prime(x_0).\;$
\end{theorem}
{\bf Proof.} Since $f$ and $g$ are intuitionistic fuzzy differentiable at $x_0$, therefore we have for any given
$\;\epsilon\;>\;0\;,\alpha\,\in\,(0,1)\;,\exists\;\delta
\,=\,\delta(\alpha,\epsilon)\;>0\;,\beta\;=\beta(\alpha,\epsilon)\;\in\;(0,1)$
such that for all $x(\neq\,x_0)\,\in\,\mathbf{R}$,
\[\mu_1(x-x_0 \;,\;\delta)\;>\;1-\beta\;\Rightarrow\;\mu_2\left(\,\frac{f(x)-f(x_0)}{x-x_0}\,-\,f^\prime (x_0) \,,\,\epsilon\,\right)\,>\,1-\alpha\]
 \[\nu_1(x-x_0 \;,\;\delta)\;<\;\beta\;\Rightarrow\;\nu_2\left(\,\frac{f(x)-f(x_0)}{x-x_0}\,-\,f^\prime (x_0) \,,\,\epsilon\,\right)\,<\,\alpha\;. \]
and
\[\mu_1(x-x_0 \;,\;\delta)\;>\;1-\beta\;\Rightarrow\;\mu_2\left(\,\frac{g(x)-g(x_0)}{x-x_0}\,-\,g^\prime (x_0) \,,\,\epsilon\,\right)\,>\,1-\alpha\]
 \[\nu_1(x-x_0 \;,\;\delta)\;<\;\beta\;\Rightarrow\;\nu_2\left(\,\frac{g(x)-g(x_0)}{x-x_0}\,-\,g ^\prime(x_0) \,,\,\epsilon\,\right)\,<\,\alpha\;. \]

Now,\[\mu_2\left(\frac{(Kf+g)(x)-(Kf+g)(x_0)}{x-x_0}-(Kf^\prime(x_0)+g^\prime(x_0)),\epsilon
\right)\hspace{6 cm}\]
\[=\mu_2\left(\frac{Kf(x)+g(x)-Kf(x_0)-g(x_0)}{x-x_0}-Kf^\prime(x_0)-g^\prime(x_0),\epsilon
\right)\hspace{6 cm}\]
\[\geq \mu_2\left(\frac{Kf(x)-Kf(x_0)}{x-x_0}-Kf^\prime(x_0),\frac{\epsilon}{2}\right)\,\ast\,
\mu_2\left(\frac{g(x)-g(x_0)}{x-x_0}-g^\prime(x_0),\frac{\epsilon}{2}\right)\hspace{6 cm}\]
\[=\mu_2\left(\frac{f(x)-f(x_0)}{x-x_0}-f^\prime(x_0),\frac{\epsilon}{2|K|}\right)\,\ast\,
\mu_2\left(\frac{g(x)-g(x_0)}{x-x_0}-g^\prime(x_0),\frac{\epsilon}{2}\right)\hspace{6 cm}\]
$>(1-\alpha)\ast(1-\alpha)=(1-\alpha),\;\;$ whenever $\,\mu_1(x-x_0 \;,\;\delta)\;>\;1-\beta\,.\\$
and \[\nu_2\left(\frac{(Kf+g)(x)-(Kf+g)(x_0)}{x-x_0}-(Kf^\prime(x_0)+g^\prime(x_0)),\epsilon\right)\hspace{6 cm}\]
\[=\nu_2\left(\frac{Kf(x)+g(x)-Kf(x_0)-g(x_0)}{x-x_0}-Kf^\prime(x_0)-g^\prime(x_0),\epsilon\right)\hspace{6 cm}\]
\[\geq \nu_2\left(\frac{Kf(x)-Kf(x_0)}{x-x_0}-Kf^\prime(x_0),\frac{\epsilon}{2}\right)\,\diamond\,
\nu_2\left(\frac{g(x)-g(x_0)}{x-x_0}-g^\prime(x_0),\frac{\epsilon}{2}\right)\hspace{6 cm}\]
\[=\nu_2\left(\frac{f(x)-f(x_0)}{x-x_0}-f^\prime(x_0),\frac{\epsilon}{2|K|}\right)\,\diamond\,
\nu_2\left(\frac{g(x)-g(x_0)}{x-x_0}-g^\prime(x_0),\frac{\epsilon}{2}\right)\hspace{6 cm}\]
$<\alpha\ast\,\alpha=\alpha,\;\;$ whenever $\,\nu_1(x-x_0 \;,\;\delta)\;<\;\beta\,.\\$
So, $Kf+g$ is intuitionistic fuzzy differentiable at $x_0\in\,\mathbf{R}$ and $(Kf+g)^\prime(x_0)\,=\,K{f}^\prime(x_0)\,+\,{g}^\prime(x_0).$

\begin{definition}
Let $\;(\;U\,,\,A\;)\;$ and $\;(\;V\,,\,B\;)\;$ be two
intuitionistic fuzzy normed linear space over the same field $\;\;k\,(=\,\mathbf{R}\,or\,\mathbf{C}).\;$ An operator $\;T\;$ from $\;(\;U\,,\,A\;)\;$ to $\;(\;V\,,\,B\;)\;$ is said to be {\bf intuitionistic fuzzy Gateaux differentiable} at $\;x_0\;\in\;U,\;$ if there exists an intuitionistic fuzzy continuous linear operator $\;G:\,\;(\;U\,,\,A\;)\longrightarrow(\;V\,,\,B\;)\;$ (generally depends upon $x_0$) and for any given
$\;\epsilon\;>\;0\;,\alpha\,\in\,(0,1)\;,\exists\;\delta
\,=\,\delta(\alpha,\epsilon)\;>0\;,\beta\;=\beta(\alpha,\epsilon)\;\in\;(0,1)$
such that for every $\;x\;\in\;U\;$ and $\;s(\neq\,0)\in\,\mathbf{R}\;$,
\[\mu_R(s,\delta)\,>\,1-\beta\;\Rightarrow\;\mu_V\left(\,\frac{T(x_0+sx)-T(x_0)}{s}\,-\,G(x) \,,\,\epsilon\,\right)\,>\,1-\alpha\]
\[\nu_R(s,\delta)\,<\,\beta\;\Rightarrow\;\nu_V\left(\,\frac{T(x_0+sx)-T(x_0)}{s}\,-\,G(x)\,,
\,\epsilon\,\right)\,<\,\alpha\;. \]
In this case, the operator $\;G\;$ is called intuitionistic fuzzy Gateaux derivative of $T$ at $x_0$ and it is denoted by $D_{f(x_0)}$.
\end{definition}

{\bf Alternative definition:} Let $\;(\;U\,,\,A\;)\;$ and $\;(\;V\,,\,B\;)\;$ be two
intuitionistic fuzzy normed linear space over the same field $\;\;k\,(=\,\mathbf{R}\,or\,\mathbf{C}).\;$ An operator $\;T\;$ from $\;(\;U\,,\,A\;)\;$ to $\;(\;V\,,\,B\;)\;$ is said to be {\bf intuitionistic fuzzy Gateaux differentiable} at $\;x_0\;\in\;U,\;$ if there exists an intuitionistic fuzzy continuous linear operator $\;G:\,\;(\;U\,,\,A\;)\longrightarrow(\;V\,,\,B\;)\;$ (generally depends upon $x_0$) such that for every $\;x\,\in\,U,\;t>0\;$ and $\;s(\neq\,0)\in\,\mathbf{R}\;$
\[\mathop {\lim }\limits_{\mu_R(s,\delta)\;\, \to\,\;1 } \mu_V\left(\,\frac{T(x_0+sx)-T(x_0)}{s}\,-\,G(x) \,,\,t\,\right)\,=\,1\]
\[\mathop {\lim }\limits_{\nu_R(s,\delta)\;\, \to\,\;0 } \nu_V\left(\,\frac{T(x_0+sx)-T(x_0)}{s}\,-\,G(x) \,,\,t\,\right)\,=\,0\]
In this case, the operator $\;G\;$ is called intuitionistic fuzzy Gateaux derivative of $T$ at $x_0$ and it is denoted by $D_{f(x_0)}$.

\begin{note}
It is easy to see that these two definitions are equivalent.
\end{note}

\begin{theorem}
Let $\;T\,:\,(U,A)\longrightarrow\,(V,B)\,$ be a linear operator, where $(U,A)$ and\, $(V,B)$ are two intuitionistic fuzzy normed linear space satisfying $(xiii)$ and $(xiv)$. If $\,T\,$ is intuitionistic fuzzy Gateaux differentiable at $\;x_0\;$ then it is unique at $\;x_0\;.$
\end{theorem}
{\bf Proof.} Let $\;G_1,\;G_2\;$ be two intuitionistic fuzzy Gateaux derivative of $T$ at $x_0$. Then for for any given
$\;\epsilon\;>\;0\;,\alpha\,\in\,(0,1)\;,\exists\;\delta
\,=\,\delta(\alpha,\epsilon)\;>0\;,\beta\;=\beta(\alpha,\epsilon)\;\in\;(0,1)$
such that for every $\;x\;\in\;U\;$ and $\;s(\neq\,0)\in\,\mathbf{R}\;$,
\[\mu_R(s,\delta)\,>\,1-\beta\;\Rightarrow\;\mu_V\left(\,\frac{T(x_0+sx)-T(x_0)}{s}\,-\,G_1(x) \,,\,\epsilon\,\right)\,>\,1-\alpha\]
\[\nu_R(s,\delta)\,<\,\beta\;\Rightarrow\;\nu_V\left(\,\frac{T(x_0+sx)-T(x_0)}{s}\,-\,G_1(x)
\,,\,\epsilon\,\right)\,<\,\alpha\;. \]
and \[\mu_R(s,\delta)\,>\,1-\beta\;\Rightarrow\;\mu_V\left(\,\frac{T(x_0+sx)-T(x_0)}{s}\,-\,G_2(x) \,,\,\epsilon\,\right)\,>\,1-\alpha\]
\[\nu_R(s,\delta)\,<\,\beta\;\Rightarrow\;
\nu_V\left(\,\frac{T(x_0+sx)-T(x_0)}{s}\,-\,G_2(x)\,,\,\epsilon\,\right)\,<\,\alpha\;. \]

\[\mu_V(\,G_1(x)-G_2(x)\,,\,t\,)\hspace{10 cm}\]
\[=\,\mu_V\left(\,\{\frac{T(x_0+sx)-T(x_0)}{s}\,-\,G_1(x)\}\,-
\,\{\frac{T(x_0+sx)-T(x_0)}{s}\,-\,G_2(x)\}\,,\,t\,\right)\]
\[\geq\;\mu_V\left(\,\frac{T(x_0+sx)-T(x_0)}{s}\,-\,G_1(x)\,,\,\frac{t}{2}\,\right)\,
\ast\,\mu_V\left(\,\frac{T(x_0+sx)-T(x_0)}{s}\,-\,G_2(x)\,,\,\frac{t}{2}\,\right)\]
\[>\,(1-\alpha)\ast(1-\alpha)\,=\,(1-\alpha)\;\;\;\;\forall\;\alpha\in(0,1).\hspace{6 cm}\]
Therefore, $\;\mu_V(\,G_1(x)-G_2(x)\,,\,t\,)\,>\,0\;\;\;\;\forall\;t>0\hspace{3 cm}(1)\\$
and \[\nu_V(\,G_1(x)-G_2(x)\,,\,t\,)\hspace{10 cm}\]
\[\leq\;\nu_V\left(\,\frac{T(x_0+sx)-T(x_0)}{s}\,-\,G_1(x)\,,\,\frac{t}{2}\,\right)
\,\diamond\,\mu_V\left(\,\frac{T(x_0+sx)-T(x_0)}{s}\,-\,G_2(x)\,,\,\frac{t}{2}\,\right)\]
\[<\,\alpha\,\diamond\,\alpha\,=\,\alpha\;\;\;\;\forall\;\alpha\in(0,1).\hspace{10 cm}\]
Therefore, $\;\nu_V(\,G_1(x)-G_2(x)\,,\,t\,)\,<\,1\;\;\;\;\forall\;t>0\hspace{3 cm}(2)\\$
From (1) and (2) we have $\;\;G_1(x)-G_2(x)\,=\theta\;.\;$
Thus $\;G_1(x)\,=\,G_2(x)\,.$

\begin{theorem}
If $\;T_1\;$ and $\;T_2\;$ have intuitionistic fuzzy Gateaux derivative at $\;x_0\;$ then $\;T\,=\,cT_1\,+\,T_2\;$ has intuitionistic fuzzy Gateaux derivative at $\;x_0,\;$ where $\;c\;$ is a scaler.
\end{theorem}
{\bf Proof.} Straight forward.

\section{Intuitionistic fuzzy Fr$\acute{e}$chet derivative}
\begin{definition}
Let $\;(\;U\,,\,A\;)\;$ and $\;(\;V\,,\,B\;)\;$ be two
intuitionistic fuzzy normed linear space over the same field $\;\;k\,(=\,\mathbf{R}\,or\,\mathbf{C}).\;$ An operator $\;T\;$ from $\;(\;U\,,\,A\;)\;$ to $\;(\;V\,,\,B\;)\;$ is said to be {\bf intuitionistic fuzzy Fr$\acute{e}$chet differentiable} at an interior $\;x_0\;\in\;U,\;$ if there exists a continuous linear operator $\;F\,:\,(\;U\,,\,A\;)\longrightarrow(\;V\,,\,B\;)$ (in general depends on $x_0$) and if for any given
$\;\epsilon\;>\;0\;,\alpha\,\in\,(0,1)\;,\exists\;\delta
\,=\,\delta(\alpha,\epsilon)\;>0\;,\beta\;=\beta(\alpha,\epsilon)\;\in\;(0,1)$
such that for all $x\;\in\;U$,
\[\mu_U(x-x_0 \;,\;\delta)\;>\;1-\beta\;\Rightarrow\;\mu_V\left(\,\frac{T(x)-T(x_0)-(x-x_0)\,F}{1-\mu_U(x-x_0 \;,\;t)}\,\,,\,\epsilon\,\right)\,>\,1-\alpha\]
\[\nu_U(x-x_0 \;,\;\delta)\;<\;\beta\;\Rightarrow\;\nu_V\left(\,\frac{T(x)-T(x_0)-(x-x_0)\,F}{\nu_U(x-x_0 \;,\;t)}\,,\,\epsilon\,\right)\,<\,\alpha\;. \]
In this case, $\;F\;$ is called intuitionistic fuzzy Fr$\acute{e}$chet derivative of $T$ at $x_0$ and is denoted by $DT(x_0)$.
\end{definition}

{\bf Alternative definition:} Let $\;(\;U\,,\,A\;)\;$ and $\;(\;V\,,\,B\;)\;$ be two
intuitionistic fuzzy normed linear space over the same field $\;k\,(=\,\mathbf{R}\,or\,\mathbf{C}).\;$ An operator $\;T\;$ from $\;(\;U\,,\,A\;)\;$ to $\;(\;V\,,\,B\;)\;$ is said to be {\bf intuitionistic fuzzy Fr$\acute{e}$chet differentiable} at an interior $\;x_0\;\in\;U,\;$ if there exists a continuous linear operator $\;F\,:\,(\;U\,,\,A\;)\longrightarrow(\;V\,,\,B\;)$ (in general depends on $x_0$) such that for every $\;t>0\;$
\[\mathop {\lim }\limits_{{\mu_U(x-x_0,t)}\;\, \to\,\;1 }\mu_V\left(\,\frac{T(x)-T(x_0)-(x-x_0)\,F}{1-\mu_U(x-x_0 \;,\;t)}\,\,,\,t\,\right)\,=\,1\]
\[\mathop {\lim }\limits_{{\nu_U(x-x_0,t)}\;\, \to\,\;0 }\mu_V\left(\,\frac{T(x)-T(x_0)-(x-x_0)\,F}{\nu_U(x-x_0 \;,\;t)}\,\,,\,t\,\right)\,=\,0\]
In this case, $\;F\;$ is called intuitionistic fuzzy Fr$\acute{e}$chet derivative of $T$ at $x_0$ and is denoted by $DT(x_0)$.

\begin{note}
It is easy to see that these two definitions are equivalent.
\end{note}

\begin{theorem}
Let $\;T\,:\,(U,A)\longrightarrow\,(V,B)\,$ be a linear operator, where $(U,A)$ and\, $(V,B)$ are two intuitionistic fuzzy normed linear space satisfying $(xiii)$ and $(xiv)$. If $\,T\,$ is intuitionistic fuzzy Fr$\acute{e}$chet differentiable at $\;x_0\;$ then it is unique at $\;x_0\;.$
\end{theorem}
{\bf Proof.} Straight forward.

\begin{example}
Let $\;U\,=\,V\,=\,\mathbf{R}\;$ and $\;[a,b]\;$ be an interval of $\;\mathbf{R}\;$ and $\;T\,:\,[a,b]\longrightarrow\mathbf{R}\,.\;$ For all $\;t>0\;$ define $\;\mu(x,t)\,=\,\frac{t}{t+\mid\,x\,\mid}\;,\;\nu(x,t)\,=\,\frac{\mid\,x\,\mid}{t+\mid\,x\,\mid}\;$ then the intuitionistic fuzzy Fr$\acute{e}$chet derivative of $\;T\;$ at $\;x_0\;$ is intuitionistic fuzzy derivative.
\end{example}
{\bf Proof.} If $\;T\;$ is intuitionistic fuzzy Fr$\acute{e}$chet differentiable at $\;x_0\;$ then for any given $\;\epsilon\;>\;0\;,\alpha\,\in\,(0,1)\;,\exists\;\delta
\,=\,\delta(\alpha,\epsilon)\;>0\;,\beta\;=\beta(\alpha,\epsilon)\;\in\;(0,1)$
such that for all $x\;\in\;U$,
\[\mu_U(x-x_0 \;,\;\delta)\;>\;1-\beta\;\Rightarrow\;\mu_V\left(\,\frac{T(x)-T(x_0)-(x-x_0)\,F}
{1-\mu_U(x-x_0 \;,\;t)}\,\,,\,\epsilon\,\right)\,>\,1-\alpha\]
\[\Rightarrow\;\mu_V\left(\,\frac{T(x)-T(x_0)-(x-x_0)\,F}{|\,x-x_0\,|}\,\,
,\,\frac{\epsilon}{t+|\,x-x_0\,|}\,\right)\,>\,1-\alpha\]
\[\Rightarrow\;\mu_V\left(\,\frac{T(x)-T(x_0)}{\,x-x_0\,}\,-\,F\,,
\,\frac{\epsilon}{t+|\,x-x_0\,|}\,\right)\,>\,1-\alpha\hspace{1.5 cm}\]
and
\[\nu_U(x-x_0 \;,\;\delta)\;<\;\beta\;\Rightarrow\;\nu_V\left(\,\frac{T(x)-T(x_0)-(x-x_0)\,F}{\nu_U(x-x_0 \;,\;t)}\,,\,\epsilon\,\right)\,<\,\alpha\;. \]
\[\Rightarrow\;\nu_V\left(\,\frac{T(x)-T(x_0)-(x-x_0)\,F}{|\,x-x_0\,|}\,\,
,\,\frac{\epsilon}{t+|\,x-x_0\,|}\,\right)\,<\,\alpha\]
\[\Rightarrow\;\nu_V\left(\,\frac{T(x)-T(x_0)}{\,x-x_0\,}\,-\,F\,,
\,\frac{\epsilon}{t+|\,x-x_0\,|}\,\right)\,<\,\alpha\hspace{1.5 cm}\]
Hence, intuitionistic fuzzy Fr$\acute{e}$chet derivative of $\;T\;$ at $\;x_0\;$ implies intuitionistic fuzzy derivative $\;T\;$ at $\;x_0\;$ and $\,T^\prime(x_0)=DT(x_0)\,.$

\begin{theorem}
An operator $\;T\;$ from $\;(\;U\,,\,A\;)\;$ to $\;(\;V\,,\,B\;)\;$ is intuitionistic fuzzy Fr$\acute{e}$chet differentiable at $\;x_0\,\in\,U\;$ then $\;T\;$ is intuitionistic fuzzy Gateaux differentiable at $\;x_0\,.$
\end{theorem}
{\bf Proof.} Since $\;T\;$ is intuitionistic fuzzy Fr$\acute{e}$chet differentiable at $\;x_0\;$, therefore we have for $\;t>0\\$
$\mu_V\left(\,\frac{T(x_0+h)-T(x_0)-DT(x_0)h}{1-\mu_U(h,t)}\,,\,t\,\right)\,>\,1-\alpha\; ,\;\;\;\;\nu_V\left(\,\frac{T(x_0+h)-T(x_0)-DT(x_0)h}{\nu_U(h,t)}\,,\,t\,\right)\,<\,\alpha\\$Now,
\[\mu_V\left(\,\frac{T(x_0+h)-T(x_0)-DT(x_0)h}{1-\mu_U(h,t)}\,,\,t\,\right)\,>\,1-\alpha \hspace{8 cm}\]
\[\Rightarrow\;\mu_V\left(\,\frac{T(x_0+sh)-T(x_0)-sDT(x_0)h}{1-\mu_U(sh,t)}\,,\,t\,\right)
\,>\,1-\alpha,\;\;\;Putting\;h=sh,\;s\neq\,0\]
\[\Rightarrow\;\mu_V\left(\,\frac{\frac{T(x_0+sh)-T(x_0)}{s}-DT(x_0)h}{\frac{1}{s}
(1-\mu_U(h,\frac{t}{|s|}))}\,,\,t\,\right)\,>\,1-\alpha\hspace{8 cm}\]
\[\Rightarrow\;\mu_V\left(\,\frac{T(x_0+sh)-T(x_0)}{s}-DT(x_0)h\,\,,\,\frac{t}{|s|}
(1-\mu_U(h,\frac{t}{|s|})\,\right)\,\,>\,1-\alpha\hspace{8 cm}\]
\[\Rightarrow\;\mu_V\left(\,\frac{T(x_0+sh)-T(x_0)}{s}-DT(x_0)h\,,\,t_1\,\right)\,\,>\,1-\alpha,\;\;\;\;where\;t_1=\,\frac{t}{|s|}(1-\mu_U(h,\frac{t}{|s|}))
\hspace{8 cm}\]
and
\[\nu_V\left(\,\frac{T(x_0+h)-T(x_0)-DT(x_0)h}{\nu_U(h,t)}\,,\,t\,\right)\,>\,1-\alpha \hspace{8 cm}\]
\[\Rightarrow\;\nu_V\left(\,\frac{T(x_0+sh)-T(x_0)-sDT(x_0)h}{\nu_U(sh,t)}\,,\,t\,\right)
\,>\,1-\alpha,\;\;\;Putting\;h=sh,\;s\neq\,0\]
\[\Rightarrow\;\nu_V\left(\,\frac{\frac{T(x_0+sh)-T(x_0)}{s}-DT(x_0)h}{\frac{1}{s}
\nu_U(h,\frac{t}{|s|}\,)}\,,\,t\,\right)\,>\,1-\alpha\hspace{8 cm}\]
\[\Rightarrow\;\nu_V\left(\,\frac{T(x_0+sh)-T(x_0)}{s}-DT(x_0)h\,\,,
\,\frac{t}{|s|}\nu_U(h,\frac{t}{|s|}\,\right)\,\,>\,1-\alpha\hspace{8 cm}\]
\[\Rightarrow\;\nu_V\left(\,\frac{T(x_0+sh)-T(x_0)}{s}-DT(x_0)h\,,\,t_2\,\right)\,\,>\,1-\alpha,\;
\;\;\;where\;t_2=\,\frac{t}{|s|}\nu_U(h,\frac{t}{|s|})
\hspace{8 cm}\]
Hence, $\;T\;$ is intuitionistic fuzzy Gateaux differentiable at $\;x_0\,$ and $\;D_{T(x_0)h}\,=\,DT(x_0)h.$

\begin{theorem}
Let $\;P:U\subset\,X\longrightarrow\,V\subset\,Y\;$ and $\;Q:V\longrightarrow\,Z\;$ be two linear operator. Suppose $\;P\;$ is intuitionistic fuzzy continuous and has intuitionistic fuzzy Gateaux derivative at $\;x_o\,\in\,U\;$ and $\;Q\;$ has intuitionistic fuzzy Fr$\acute{e}$chet derivative at $\;y_0\,=\,P(x_0)$. Then $\;R\,=\,QP\;$ has intuitionistic fuzzy Gateaux derivative at $\;x_0\;$ and $\;D_{R(x_0)}\,=\,DQ(y_0)\,D_{P(x_0)}$.
\end{theorem}
{\bf Proof.} We write $\;G\,=\,D_{P(x_0)}\;$ and $\;F\,=\,DQ(y_0)\;$ for shortness. Let $\;x\,\in\,X\;$ and we further write $\;\triangle y\,=\,P(x_0+sx)-P(x_0).\;$ Then
\[\mu\left(\frac{\,R(x_0+sx)-R(x_0)}{s}-FG\,,\,t\right)\,=
\,\mu\left(\frac{\,QP(x_0+sx)-QP(x_0)}{s}-FG\,,\,t\right)\,\]
\[=\mu\left(\frac{\,F(\triangle y)+A(\triangle y)}{s}-FG\,,\,t\right)\,\;\;,where\;A(\triangle y)\,=\,Q(y_0+\triangle y)-Q(y_0)-F(\triangle y)\]
\[=\mu\left(F \frac{P(x_0+sx)-P(x_0)}{s}+\frac{A(\triangle y)}{s}-FG\,,\,t\right)\hspace{7 cm}\]
\[\geq\,\mu\left(F \frac{P(x_0+sx)-P(x_0)}{s}-FG\,,\,\frac{t}{2}\right)\,\ast\,\mu\left(\frac{A(\triangle y)}{\mu(\triangle y,t_1)}\,\frac{\mu(P(x_0+sx)-P(x_0),t_1)}{s}\,,\,\frac{t}{2}\right)\]
\[=\mu\left(\frac{P(x_0+sx)-P(x_0)}{s}-G\,,\,\frac{t}{2\mu(F,t_2)}\right)\,\ast\,
\mu\left(\frac{A(\triangle y)}{\mu(\triangle y,t_1)}\,,\,\frac{ts}{2\mu(P(x_0+sx)-P(x_0),t_1)}\right)\]
\[=\mu\left(\frac{P(x_0+sx)-P(x_0)}{s}-G\,,\,\frac{t}{2\mu(F,t_2)}\right)\,\hspace{8 cm}\]\[\ast\,\mu\left(\frac{Q(y_0+\triangle y)-Q(y_0)-F(\triangle y)}{\mu(\triangle y,t_1)}\,,\,\frac{ts}{2\mu(P(x_0+sx)-P(x_0),t_1)}\right)\hspace{-3 cm}\]
\[>\,(1-\alpha)\,\ast\,(1-\alpha)\;=\,(1-\alpha)\hspace{10 cm}\] since $\,P\,$ has intuitionistic fuzzy Gateaux derivative and $\,Q\,$ has intuitionistic fuzzy Fr$\acute{e}$chet derivative.\\
and \[\nu\left(\frac{\,R(x_0+sx)-R(x_0)}{s}-FG\,,\,t\right)\,=\,
\nu\left(\frac{\,QP(x_0+sx)-QP(x_0)}{s}-FG\,,\,t\right)\,\]
\[=\nu\left(\frac{\,F(\triangle y)+A(\triangle y)}{s}-FG\,,\,t\right)\,\;\;,where\;A(\triangle y)\,=\,Q(y_0+\triangle y)-Q(y_0)-F(\triangle y)\]
\[=\nu\left(F \frac{P(x_0+sx)-P(x_0)}{s}+\frac{A(\triangle y)}{s}-FG\,,\,t\right)\hspace{7 cm}\]
\[\leq\,\nu\left(F \frac{P(x_0+sx)-P(x_0)}{s}-FG\,,\,\frac{t}{2}\right)\,\diamond\,\nu\left(\frac{A(\triangle y)}{1-\nu(\triangle y,t_1)}\,\frac{1-\nu(P(x_0+sx)-P(x_0),t_1)}{s}\, ,\,\frac{t}{2}\right)\]
\[=\nu\left(\frac{P(x_0+sx)-P(x_0)}{s}-G\,,\,\frac{t}{2\nu(F,t_2)}\right)\,\hspace{10 cm}\]\[\diamond\,
\nu\left(\frac{A(\triangle y)}{1-\nu(\triangle y,t_1)}\, ,\,\frac{ts}{2(1-\nu(P(x_0+sx)-P(x_0),t_1))}\right)\]
\[=\nu\left(\frac{P(x_0+sx)-P(x_0)}{s}-G\,,\,\frac{t}{2\nu(F,t_2)}\right)\,\hspace{8 cm}\]\[\diamond\,\nu\left(\frac{Q(y_0+\triangle y)-Q(y_0)- F(\triangle y)}{1-\nu(\triangle y,t_1)}\,,\,\frac{ts}{2(1-\nu(P(x_0+sx)-P(x_0),t_1))}\right)\hspace{-3 cm}\]
\[<\,\alpha\,\diamond\,\alpha\;=\,\alpha\,.\hspace{15 cm}\]
Since $\,P\,$ has intuitionistic fuzzy Gateaux derivative and $\,Q\,$ has intuitionistic fuzzy Fr$\acute{e}$chet derivative.
Hence $\;R\,=\,QP\;$ has intuitionistic fuzzy Gateaux derivative at $\;x_0\;$ and $\;D_{R(x_0)}\,=\,DQ(y_0)\,D_{P(x_0)}$.

\end{document}